\documentclass[12pt]{article}
\usepackage{amssymb}
\usepackage{amsfonts}
\usepackage{amsmath}

\begin{document}

\begin{center}
{\Large Real matrix representations for the complex quaternions }%
\begin{equation*}
\end{equation*}%
\begin{equation*}
\end{equation*}

Cristina FLAUT \ and $\ $Vitalii \ SHPAKIVSKYI%
\begin{equation*}
\end{equation*}%
\begin{equation*}
\end{equation*}
\end{center}

\textbf{Abstract.} {\small Starting from known results, due to Y. Tian in
[Ti; 00], referring to the real matrix representations of the real
quaternions,}\ {\small in this paper \ we will \ investigate the \ left and
\ right real matrix representations for \ the complex quaternions and we
will give some examples in the special case of the complex Fibonacci
quaternions.}%
\begin{equation*}
\end{equation*}



\begin{equation*}
\end{equation*}

\textbf{1. Introduction}%
\begin{equation*}
\end{equation*}

We know that each finite dimensional associative \ algebra $A$ over an
arbitrary field $K$ is isomorphic with a subalgebra of the algebra $\mathcal{%
M}_{n}\left( K\right) $, with $n=\dim _{K}A$. Therefore, we can find a
faithful representation of the algebra $A$ in the algebra $\mathcal{M}%
_{n}\left( K\right) .$ For example, the real quaternion division algebra is
algebraically isomorphic to a $4\times 4$ real matrix algebra. Starting from
some results obtained by Y. Tian in [Ti; 00] and in [Ti; 00(1)]{\small , }in
this paper we will show that the complex quaternion algebra is algebraically
isomorphic to a $8\times 8$ real matrix algebra and will investigate the
properties of the obtained left and right real matrix representations for \
the complex quaternions. In Section 3, we will provide some examples in the
special case of the complex Fibonacci quaternions.

Let $K$ be the field $\{\left( 
\begin{array}{cc}
a & -b \\ 
b & a%
\end{array}%
\right) $ $\mid $ $a,b\in \mathbb{R}\}.$ The map 
\begin{equation*}
\varphi :\mathbb{C}\rightarrow K,\varphi \left( a+bi\right) =\left( 
\begin{array}{cc}
a & -b \\ 
b & a%
\end{array}%
\right) ,
\end{equation*}%
where $i^{2}=-1$ is a fields morphism and $\varphi \left( z\right) =\left( 
\begin{array}{cc}
a & -b \\ 
b & a%
\end{array}%
\right) $ is called the matrix representation of the element $z=a+bi\in 
\mathbb{C}.$

Let $\mathbb{H}$ be the real\ division quaternion algebra, the algebra of
the elements of the form $a=a_{0}+a_{1}i+a_{2}j+a_{3}k,$ where

\begin{equation*}
a_{i}\in \mathbb{R},i^{2}=j^{2}=k^{2}=-1,
\end{equation*}%
and \ 
\begin{equation*}
ij=-ji=k,jk=-kj=i,ki=-ik=j.
\end{equation*}%
$\mathbb{H}\ $is an algebra $\ $over the field $\mathbb{R}.$ The set $%
\{1,i,j,k\}$ is a basis in $\mathbb{H}$. The conjugate of the real
quaternion $a=a_{0}+a_{1}i+a_{2}j+a_{3}k\,$\ is the quaternion $\overline{a}%
=a_{0}-a_{1}i-a_{2}j-a_{3}k$ and $n\left( a\right) =a\overline{a}=\overline{a%
}a$ is called \textit{the norm }of the real quaternion $a.$

A complex quaternion is an element of the form $%
Q=c_{0}+c_{1}e_{1}+c_{2}e_{2}+c_{3}e_{3},$ where $c_{n}\in \mathbb{C},n\in
\{0,1,2,3\},$%
\begin{equation*}
\ e_{n}^{2}=-1,\,\,\,n\in \{1,2,3\}
\end{equation*}%
and \ 
\begin{equation*}
e_{m}e_{n}=-e_{n}e_{m}=\beta _{mn}e_{t},\,\,\beta _{mn}\in \{-1,1\},m\neq
n,m,n\in \{1,2,3\},
\end{equation*}%
$\beta _{mn}$ and $e_{t}$ being uniquely determined by $e_{m}$ and $e_{n}.$
We denote by $\mathbb{H}_{C}$ the algebra of the complex quaternions\textit{,%
} called \textit{the complex quaternion algebra}. This algebra is an algebra
over the field $\mathbb{C}.$ The set $\{1,e_{1},e_{2},e_{3}\}$ is a basis in 
$\mathbb{H}_{C}$.

The map $\gamma :\mathbb{R}\rightarrow \mathbb{C},\gamma \left( a\right) =a$
is the inclusion morphism between $\mathbb{R}$-algebras $\mathbb{R}$ and $%
\mathbb{C}.$ We denote by $\mathbb{F}$ the $\mathbb{C}$-subalgebra of the
algebra $\mathbb{H}_{C},$%
\begin{equation*}
\mathbb{F}=\{Q\in \mathbb{H}_{C}\ \mid
~Q=c_{0}+c_{1}e_{1}+c_{2}e_{2}+c_{3}e_{3},c_{n}\in \mathbb{R},n\in
\{0,1,2,3\}\}.
\end{equation*}

\bigskip By the scalar restriction, $\mathbb{F}$ became an algebra over $%
\mathbb{R},$ with the multiplication $"\cdot "$ 
\begin{equation*}
a\cdot Q=\gamma \left( a\right) Q=aQ,a\in \mathbb{R},Q\in \mathbb{F}.
\end{equation*}%
We denote this algebra by $\mathbb{H}_{R}.$ The map 
\begin{equation*}
\delta :\mathbb{H\rightarrow H}_{R},\delta \left( 1\right) =1,\delta \left(
i\right) =e_{1},\delta \left( j\right) =e_{2},\delta \left( k\right) =e_{3}
\end{equation*}%
and%
\begin{equation*}
\delta \left( a_{0}+a_{1}i+a_{2}j+a_{3}k\right)
=a_{0}+a_{1}e_{1}+a_{2}e_{2}+a_{3}e_{3},
\end{equation*}%
where $a_{m}\in \mathbb{R},m\in \{0,1,2,3\}$ is an algebra isomorphism
between the algebras $\mathbb{H}$ and $\mathbb{H}_{R}.$The algebra $\mathbb{H%
}_{R}$ has the same basis $\{1,e_{1},e_{2},e_{3}\}$ as the algebra $\mathbb{H%
}_{C}.$ From now one, we will identify the quaternion $%
a_{0}+a_{1}i+a_{2}j+a_{3}k$ with the "complex" quaternion $%
a_{0}+a_{1}e_{1}+a_{2}e_{2}+a_{3}e_{3},$ $a_{m}\in \mathbb{R},m\in
\{0,1,2,3\}$ and instead of $\mathbb{H}_{R}$ we will use $\mathbb{H}.$

It results that the element $Q\in \mathbb{H}%
_{C},~Q=c_{0}+c_{1}e_{1}+c_{2}e_{2}+c_{3}e_{3},c_{m}\in \mathbb{C},m\in
\{0,1,2,3\},$ can be written as \ $%
Q=(a_{0}+ib_{0})+(a_{1}+ib_{1})e_{1}+(a_{2}+ib_{2})e_{2}+(a_{3}+ib_{3})e_{3}, 
$ where $a_{m},b_{m}\in \mathbb{R},m\in \{0,1,2,3\}$ and $i^{2}=-1.$

Therefore, we can write a complex quaternion under the form 
\begin{equation*}
Q=a+ib,
\end{equation*}%
with $a,b\in \mathbb{H}$, $%
a=a_{0}+a_{1}e_{1}+a_{2}e_{2}+a_{3}e_{3},b=b_{0}+b_{1}e_{1}+b_{2}e_{2}+b_{3}e_{3}. 
$

The conjugate of the complex quaternion $Q$ is the element $\overline{Q}%
=c_{0}-c_{1}e_{1}-c_{2}e_{2}-c_{3}e_{3}$. It results that 
\begin{equation}
\overline{Q}=\overline{a}+i\overline{b}.  \tag{1.1.}
\end{equation}

For the quaternion $a=a_{0}+a_{1}e_{1}+a_{2}e_{2}+a_{3}e_{3}\in \mathbb{H}$,
we define the element%
\begin{equation}
a^{\ast }=a_{0}+a_{1}e_{1}-a_{2}e_{2}-a_{3}e_{3}.  \tag{1.2.}
\end{equation}%
We remark that 
\begin{equation}
(a^{\ast })^{\ast }=a  \tag{1.3.}
\end{equation}%
and%
\begin{equation}
\left( a+b\right) ^{\ast }=a^{\ast }+b^{\ast },  \tag{1.4.}
\end{equation}%
for all $a,b\in \mathbb{H}.$

For \ the quaternion algebra $\mathbb{H},\,$in [Ti; 00], was defined $\ $the
map \ 
\begin{equation*}
\lambda :\mathbb{H}\rightarrow \mathcal{M}_{4}\left( \mathbb{R}\right)
,\lambda \left( a\right) =\left( 
\begin{array}{llll}
a_{0} & -a_{1} & -a_{2} & -a_{3} \\ 
a_{1} & a_{0} & -a_{3} & a_{2} \\ 
a_{2} & a_{3} & a_{0} & -a_{1} \\ 
a_{3} & -a_{2} & a_{1} & a_{0}%
\end{array}%
\right) ,
\end{equation*}%
where $a=a_{0}+a_{1}e_{1}+a_{2}e_{2}+a_{3}e_{3}\in \mathbb{H},\,\,$\ is an
isomorphism between $\mathbb{H}\,\ $\ and the algebra of the \ matrices: 
\begin{equation*}
\left\{ \left( 
\begin{array}{llll}
a_{0} & -a_{1} & -a_{2} & -a_{3} \\ 
a_{1} & a_{0} & -a_{3} & a_{2} \\ 
a_{2} & a_{3} & a_{0} & -a_{1} \\ 
a_{3} & -a_{2} & a_{1} & a_{0}%
\end{array}%
\right) ,a_{0},a_{1},a_{2},a_{3}\in \mathbb{R}\right\} .\,\,
\end{equation*}

We remark that the matrix $\lambda \left( a\right) \in \mathcal{M}_{4}\left( 
\mathbb{R}\right) $ has as columns the coefficients in $\mathbb{R}$ \ of the
basis $\{1,e_{1},e_{2},e_{3}\}$ for the elements $\
\{a,ae_{1},ae_{2},ae_{3}\}.$

The matrix $\lambda \left( a\right) $ is called \textit{the left matrix
representation} of the element $a\in \mathbb{H}.\medskip $

Analogously with the left matrix representation, \ for the element $a\in 
\mathbb{H}\,,$ in [Ti; 00], was defined$\ \ $\ \textit{the right matrix
representation}:

\begin{equation*}
\rho :\mathbb{H}\rightarrow \mathcal{M}_{4}\left( \mathbb{R}\right)
,\,\,\,\rho \left( a\right) =\left( 
\begin{array}{llll}
a_{0} & -a_{1} & -a_{2} & -a_{3} \\ 
a_{1} & a_{0} & a_{3} & -a_{2} \\ 
a_{2} & -a_{3} & a_{0} & a_{1} \\ 
a_{3} & a_{2} & -a_{1} & a_{0}%
\end{array}%
\right) ,\medskip
\end{equation*}%
where $a=a_{0}+a_{1}e_{1}+a_{2}e_{2}+a_{3}e_{3}\in $ $\mathbb{H}.\medskip $ $%
\,$

We remark that the matrix $\rho \left( a\right) \in \mathcal{M}_{4}\left( 
\mathbb{R}\right) $ has as columns the coefficients in $\mathbb{R}$ \ of the
basis $\{1,e_{1},e_{2},e_{3}\}$ for the elements $\
\{a,e_{1}a,e_{2}a,e_{3}a\}.$\medskip\ 

\textbf{Proposition 1.1. }[Ti; 00] \textit{For }$x,y\in \mathbb{H}$\textit{\
and }$r\in K$ \textit{\ we have:}

\textit{i) }$\lambda \left( x+y\right) =\lambda \left( x\right) +\lambda
\left( y\right) ,\,\lambda \left( xy\right) =\lambda \left( x\right) \lambda
\left( y\right) ,\,\lambda \left( rx\right) =r\lambda \left( x\right) ,$

$\lambda \left( 1\right) =I_{4},\,r\in K.$

\textit{ii) }$\rho \left( x+y\right) =\rho \left( x\right) +\rho \left(
y\right) ,\rho \left( xy\right) =\rho \left( y\right) \rho \left( x\right)
,\rho \left( rx\right) =r\rho \left( x\right) ,$

$\rho \left( 1\right) =I_{4},\,r\in K.$

\textit{iii) }$\lambda \left( x^{-1}\right) =\left( \lambda \left( x\right)
\right) ^{-1},\rho \left( x^{-1}\right) =\left( \rho \left( x\right) \right)
^{-1},\,\,$\textit{for }$x\neq 0.\Box \medskip $

\textbf{Proposition 1.2.} [Ti; 00] \textit{For }$x\in \mathbb{H},$ \textit{%
let} $\overrightarrow{x}=\left( a_{0},a_{1},a_{2},a_{3}\right) ^{t}$\textit{%
\ }$\in \mathcal{M}_{1\times 4}\left( K\right) ,$\textit{\ be the vector
representation of the element }$x$.\textit{\ Therefore for all }$a,b,x\in 
\mathbb{H}$\textit{\ the following relations are fulfilled:}

\textit{i) }$\overrightarrow{ax}=\lambda \left( a\right) \overrightarrow{x}.$

\textit{ii) }$\overrightarrow{xb}=\rho \left( b\right) \overrightarrow{x}.$

\textit{iii) }$\overrightarrow{axb}=\lambda \left( a\right) \rho \left(
b\right) \overrightarrow{x}=\rho \left( b\right) \lambda \left( a\right) 
\overrightarrow{x}.$

\textit{iv)}$\rho \left( b\right) \lambda \left( a\right) =\lambda \left(
a\right) \rho \left( b\right) .$

\textit{v)} $\det \left( \lambda \left( x\right) \right) =\det \left( \rho
\left( x\right) \right) =\left( n\left( x\right) \right) ^{2}.\medskip \Box
\medskip $

For details about the matrix representations of the real quaternions, the
reader is referred to [Ti; 00].

\medskip\ 

\begin{equation*}
\end{equation*}
\textbf{2. Main results}%
\begin{equation*}
\end{equation*}

Let $\theta $ \ be the matrix $\theta =\left( 
\begin{array}{cccc}
0 & -1 & 0 & 0 \\ 
1 & 0 & 0 & 0 \\ 
0 & 0 & 0 & -1 \\ 
0 & 0 & 1 & 0%
\end{array}%
\right) =\lambda \left( e_{1}\right) =\lambda \left( i\right) .$ The matrix 
\begin{equation*}
\Gamma \left( Q\right) =\left( 
\begin{array}{cc}
\lambda \left( a\right) & -\lambda \left( b^{\ast }\right) \\ 
\lambda \left( b\right) & \lambda \left( a^{\ast }\right)%
\end{array}%
\right) ,
\end{equation*}%
where $Q=a+ib$ is a complex quaternion, with $%
a=a_{0}+a_{1}e_{1}+a_{2}e_{2}+a_{3}e_{3}\in \mathbb{H}%
,b=b_{0}+b_{1}e_{1}+b_{2}e_{2}+b_{3}e_{3}\in \mathbb{H}$ and $i^{2}=-1,$ is
called \textit{the left real matrix representation for the \ complex
quaternion} $Q.$ The \textit{right real matrix representation} for the \
complex quaternion $Q$ is the matrix: 
\begin{equation*}
\Theta \left( Q\right) =\left( 
\begin{array}{cc}
\rho \left( a\right) & -\rho \left( b\right) \\ 
\rho \left( b^{\ast }\right) & \rho \left( a^{\ast }\right)%
\end{array}%
\right) .
\end{equation*}%
We remark that \ $\Gamma \left( Q\right) ,\Theta \left( Q\right) \in 
\mathcal{M}_{8}\left( \mathbb{R}\right) .$

Now, let $M$ be the matrix 
\begin{equation*}
M=\left( 1,-e_{1},-e_{2},-e_{3}\right) ^{t}.
\end{equation*}%
$\medskip $

\textbf{Proposition 2.1.} \textit{If \ }$%
a=a_{0}+a_{1}e_{1}+a_{2}e_{2}+a_{3}e_{3}\in \mathbb{H},$ \textit{we have:}

\textit{i)} \ $\lambda \left( a\right) M=Ma.$

\textit{ii)} $\theta M=Me_{1}.$

\textit{iii)} $\lambda \left( ia\right) =\theta \lambda \left( a\right) $ 
\textit{and} $\ \lambda \left( ai\right) =\lambda \left( a\right) \theta
.\medskip $

\textbf{Proof.} i) \ $\lambda \left( a\right) M$=$\left( 
\begin{array}{llll}
a_{0} & -a_{1} & -a_{2} & -a_{3} \\ 
a_{1} & a_{0} & -a_{3} & a_{2} \\ 
a_{2} & a_{3} & a_{0} & -a_{1} \\ 
a_{3} & -a_{2} & a_{1} & a_{0}%
\end{array}%
\right) \left( 
\begin{array}{c}
1 \\ 
-e_{1} \\ 
-e_{2} \\ 
-e_{3}%
\end{array}%
\right) $=\newline
=$\left( 
\begin{array}{c}
a_{0}\text{+}a_{1}e_{1}\text{+}a_{2}e_{2}\text{+}a_{3}e_{3} \\ 
a_{1}\text{-}a_{0}e_{1}\text{+}a_{3}e_{2}\text{-}a_{2}e_{3} \\ 
a_{2}\text{-}a_{3}e_{1}\text{-}a_{0}e_{2}\text{+}a_{1}e_{3} \\ 
a_{3}\text{+}a_{2}e_{1}\text{-}a_{1}e_{2}\text{-}a_{0}e_{3}%
\end{array}%
\right) $=$\left( 
\begin{array}{c}
a_{0}\text{+}a_{1}e_{1}\text{+}a_{2}e_{2}\text{+}a_{3}e_{3} \\ 
\text{-}e_{1}(a_{0}\text{+}a_{1}e_{1}\text{+}a_{2}e_{2}\text{+}a_{3}e_{3})
\\ 
\text{-}e_{2}\left( a_{0}\text{+}a_{1}e_{1}\text{+}a_{2}e_{2}\text{+}%
a_{3}e_{3}\right) \\ 
\text{-}e_{3}\left( a_{0}\text{+}a_{1}e_{1}\text{+}a_{2}e_{2}\text{+}%
a_{3}e_{3}\right)%
\end{array}%
\right) $=\newline
=$\left( 
\begin{array}{c}
1 \\ 
-e_{1} \\ 
-e_{2} \\ 
-e_{3}%
\end{array}%
\right) a=Ma.$

ii) $\theta M=\left( 
\begin{array}{cccc}
0 & -1 & 0 & 0 \\ 
1 & 0 & 0 & 0 \\ 
0 & 0 & 0 & -1 \\ 
0 & 0 & 1 & 0%
\end{array}%
\right) \left( 
\begin{array}{c}
1 \\ 
-e_{1} \\ 
-e_{2} \\ 
-e_{3}%
\end{array}%
\right) $=$\left( 
\begin{array}{c}
e_{1} \\ 
1 \\ 
e_{3} \\ 
-e_{2}%
\end{array}%
\right) $=$\newline
$=$\left( 
\begin{array}{c}
1 \\ 
-e_{1} \\ 
-e_{2} \\ 
-e_{3}%
\end{array}%
\right) e_{1}=Me_{1}.$

iii) \ For $a=a_{0}+a_{1}e_{1}+a_{2}e_{2}+a_{3}e_{3}\in \mathbb{H},$ we have 
$ia=-a_{1}+a_{0}e_{1}-a_{3}e_{2}+a_{2}e_{3}.$ It result that 
\begin{equation*}
\lambda \left( ia\right) \text{=}\left( 
\begin{array}{llll}
-a_{1} & -a_{0} & a_{3} & -a_{2} \\ 
a_{0} & -a_{1} & -a_{2} & -a_{3} \\ 
-a_{3} & a_{2} & -a_{1} & -a_{0} \\ 
a_{2} & a_{3} & a_{0} & -a_{1}%
\end{array}%
\right) .
\end{equation*}%
Since $\theta \lambda \left( a\right) $=\ $\left( 
\begin{array}{cccc}
0 & \text{-}1 & 0 & 0 \\ 
1 & 0 & 0 & 0 \\ 
0 & 0 & 0 & \text{-}1 \\ 
0 & 0 & 1 & 0%
\end{array}%
\right) $ $\left( 
\begin{array}{llll}
a_{0} & -a_{1} & -a_{2} & -a_{3} \\ 
a_{1} & a_{0} & -a_{3} & a_{2} \\ 
a_{2} & a_{3} & a_{0} & -a_{1} \\ 
a_{3} & -a_{2} & a_{1} & a_{0}%
\end{array}%
\right) $=\newline
=$\left( 
\begin{array}{llll}
-a_{1} & -a_{0} & a_{3} & -a_{2} \\ 
a_{0} & -a_{1} & -a_{2} & -a_{3} \\ 
-a_{3} & a_{2} & -a_{1} & -a_{0} \\ 
a_{2} & a_{3} & a_{0} & -a_{1}%
\end{array}%
\right) ,$ we obtain the asked relation.$~\Box $

\medskip \medskip

\textbf{Proposition 2.2.} \textit{Let} $a,x\in \mathbb{H}$ \textit{be two
quaternions, then the following relations are true:}

\textit{i)} $a^{\ast }i=ia,$ \textit{where} $i^{2}=-1.$

\textit{ii)} $ai=ia^{\ast },\ $\ \textit{where} $i^{2}=-1.$

\textit{iii)} $-a^{\ast }=iai,$ \textit{where} $i^{2}=-1.$

\textit{iv)} $\left( xa\right) ^{\ast }=x^{\ast }a^{\ast }.$

\textit{v)} \textit{For} $X,A\in \mathbb{H}_{C},X=x+iy,A=a+ib,$ \textit{we
have} 
\begin{equation*}
XA=xa-y^{\ast }b+i\left( x^{\ast }b+ya\right) .\medskip
\end{equation*}

\textbf{Proof}. \ Relations from i), ii), iii) are obviously.

iv) From ii), it results $\left( xa\right) ^{\ast }=-i\left( xa\right)
i=-ixai=(ixi)(iai)=x^{\ast }a^{\ast }.$

v) We obtain\newline
$XA=\left( x+iy\right) \left( a+ib\right) =xa+xib+iya+iyib=$\newline
$=xa-y^{\ast }b+i\left( x^{\ast }b+ya\right) .\Box \medskip $

\textbf{Proposition 2.3.} \textit{For} $X,A\in \mathbb{H}_{C},X=x+iy,A=a+ib,$
\textit{we have }$\Gamma \left( XA\right) =\Gamma \left( X\right) \Gamma
\left( A\right) .\medskip $

\textbf{Proof.} From Proposition 1.2 i) \ and Proposition 2.2  iv), it
results that\newline
$\Gamma \left( X\right) \Gamma \left( A\right) $=$\left( 
\begin{array}{cc}
\lambda \left( x\right)  & -\lambda \left( y^{\ast }\right)  \\ 
\lambda \left( y\right)  & \lambda \left( x^{\ast }\right) 
\end{array}%
\right) \left( 
\begin{array}{cc}
\lambda \left( a\right)  & -\lambda \left( b^{\ast }\right)  \\ 
\lambda \left( b\right)  & \lambda \left( a^{\ast }\right) 
\end{array}%
\right) $=\medskip \newline
=$\left( 
\begin{array}{cc}
\lambda \left( x\right) \lambda \left( a\right) -\lambda \left( y^{\ast
}\right) \lambda \left( b\right)  & \text{-}\lambda \left( x\right) \lambda
\left( b^{\ast }\right) \text{-}\lambda \left( y^{\ast }\right) \lambda
\left( a^{\ast }\right)  \\ 
\lambda \left( y\right) \lambda \left( a\right) +\lambda \left( x^{\ast
}\right) \lambda \left( b\right)  & \text{-}\lambda \left( y\right) \lambda
\left( b^{\ast }\right) \text{+}\lambda \left( x^{\ast }\right) \lambda
\left( a^{\ast }\right) 
\end{array}%
\right) $=\medskip \newline
=$\left( 
\begin{array}{cc}
\lambda (xa-y^{\ast }b) & -\lambda (xb^{\ast }+y^{\ast }a^{\ast }) \\ 
\lambda (ya+x^{\ast }b) & \lambda (-yb^{\ast }+x^{\ast }a^{\ast })%
\end{array}%
\right) .\medskip $

$\Gamma \left( XA\right) $=$\left( 
\begin{array}{cc}
\lambda (xa-y^{\ast }b) & -\lambda ((x^{\ast }b+ya)^{\ast }) \\ 
\lambda (x^{\ast }b+ya) & \lambda ((xa-y^{\ast }b)^{\ast })%
\end{array}%
\right) $=\medskip \newline
=$\left( 
\begin{array}{cc}
\lambda (xa-y^{\ast }b) & -\lambda (xb^{\ast }+y^{\ast }a^{\ast }) \\ 
\lambda (ya+x^{\ast }b) & \lambda (x^{\ast }a^{\ast }-yb^{\ast })%
\end{array}%
\right) $.$~\Box \medskip \medskip $

\textbf{Definition 2.4}. For $X\in \mathbb{H}_{C},X=x+iy,$ we denote by 
\begin{equation*}
\overrightarrow{X}=(\overrightarrow{x},\overrightarrow{y})^{t}\in \mathcal{M}%
_{8\times 1}\left( \mathbb{R}\right)
\end{equation*}%
$~$\newline
\textit{the vector representation }of the element $X$, where\newline
$x$=$x_{0}$+$x_{1}e_{1}$+$x_{2}e_{2}$+$x_{3}e_{3}\in \mathbb{H},y$=$y_{0}$+$%
y_{1}e_{1}$+$y_{2}e_{2}$+$y_{3}e_{3}\in \mathbb{H}$ and\newline
$\overrightarrow{x}$=$(x_{0},x_{1},x_{2},x_{3})^{t}\in \mathcal{M}_{4\times
1}\left( \mathbb{R}\right) ,$\newline
$\overrightarrow{y}$=$(y_{0},y_{1},y_{2},y_{3})^{t}\in \mathcal{M}_{4\times
1}\left( \mathbb{R}\right) $ are the vector representations for the
quaternions $x$ and $y,$ as was defined in Proposition 1.2.\medskip

\textbf{Proposition 2.5.} \textit{Let} $X\in \mathbb{H}_{C},X=x+iy,x,y\in 
\mathbb{H},$ \textit{then:}

\textit{i)} $\overrightarrow{X}=\Gamma \left( X\right) \left( 
\begin{array}{c}
1 \\ 
0%
\end{array}%
\right) ,$ \textit{where} $1=I_{4}\in \mathcal{M}_{4}\left( \mathbb{R}%
\right) $ \textit{is the identity matrix and } $0=O_{4}\in \mathcal{M}%
_{4}\left( \mathbb{R}\right) $ \textit{is the zero matrix.}

\textit{ii)} $\overrightarrow{AX}=\Gamma \left( A\right) \overrightarrow{X}.$

\textit{iii)} $\alpha \overrightarrow{y^{\ast }}=\overrightarrow{y},$ 
\textit{where} $\alpha $=$\left( 
\begin{array}{cccc}
1 & 0 & 0 & 0 \\ 
0 & 1 & 0 & 0 \\ 
0 & 0 & \text{-}1 & 0 \\ 
0 & 0 & 0 & \text{-}1%
\end{array}%
\right) \in \mathcal{M}_{4}\left( \mathbb{R}\right) .$

\textit{iv)} $\alpha ^{2}=I_{4}.\medskip $

\textbf{Proof.} i) $\Gamma \left( X\right) \left( 
\begin{array}{c}
1 \\ 
0%
\end{array}%
\right) $=$\left( 
\begin{array}{cc}
\lambda \left( x\right) & \text{-}\lambda \left( y^{\ast }\right) \\ 
\lambda \left( y\right) & \lambda \left( x^{\ast }\right)%
\end{array}%
\right) \left( 
\begin{array}{c}
1 \\ 
0%
\end{array}%
\right) $=$\left( 
\begin{array}{c}
\lambda \left( x\right) \\ 
\lambda \left( y\right)%
\end{array}%
\right) $=\medskip \newline
=$\left( 
\begin{array}{c}
\lambda \left( 1\cdot x\right) \\ 
\lambda \left( 1\cdot y\right)%
\end{array}%
\right) $=$\left( 
\begin{array}{c}
\lambda \left( 1\right) \overrightarrow{x} \\ 
\lambda \left( 1\right) \overrightarrow{y}%
\end{array}%
\right) $=$\left( 
\begin{array}{c}
\overrightarrow{x} \\ 
\overrightarrow{y}%
\end{array}%
\right) .\medskip $

ii) From i), we obtain that\newline
$\overrightarrow{AX}$=$\Gamma \left( AX\right) \left( 
\begin{array}{c}
1 \\ 
0%
\end{array}%
\right) $=$\Gamma \left( A\right) \Gamma \left( X\right) \left( 
\begin{array}{c}
1 \\ 
0%
\end{array}%
\right) $=$\Gamma \left( A\right) \overrightarrow{X}.\medskip $

iii) \ $\alpha \overrightarrow{y^{\ast }}$=$\left( 
\begin{array}{cccc}
1 & 0 & 0 & 0 \\ 
0 & 1 & 0 & 0 \\ 
0 & 0 & -1 & 0 \\ 
0 & 0 & 0 & -1%
\end{array}%
\right) \left( 
\begin{array}{c}
y_{0} \\ 
y_{1} \\ 
-y_{2} \\ 
-y_{3}%
\end{array}%
\right) $= $\left( 
\begin{array}{c}
y_{0} \\ 
y_{1} \\ 
y_{2} \\ 
y_{3}%
\end{array}%
\right) $=$\overrightarrow{y}.\Box \medskip $

\textbf{Proposition 2.6.} \textit{Let} $M_{8}$ \ \textit{be the matrix} $%
M_{8}=\left( 
\begin{array}{c}
\theta M \\ 
-M%
\end{array}%
\right) $\textit{, therefore we have} $-\frac{1}{4}M_{8}^{t}M_{8}=1.\medskip 
$

\textbf{Proof.} It results \newline
$M_{8}^{t}M_{8}$=$\left( 
\begin{array}{cccccccc}
e_{1} & -1 & e_{3} & e_{2} & -1 & e_{1} & e_{2} & e_{3}%
\end{array}%
\right) \left( 
\begin{array}{c}
e_{1} \\ 
-1 \\ 
e_{3} \\ 
e_{2} \\ 
-1 \\ 
e_{1} \\ 
e_{2} \\ 
e_{3}%
\end{array}%
\right) =-4.\Box $

\bigskip

\textbf{Theorem 2.7}. \textit{Let }$Q\in \mathbb{H}_{C}$ \textit{be a
complex quaternion.} \textit{With the above notations, the following
relations are fulfilled:}

\textit{i)} $\Gamma ^{t}\left( Q^{\ast }\right) M_{8}=M_{8}Q,$ \textit{where}
$Q=x+iy,Q^{\ast }=x^{\ast }+iy,x,y\in \mathbb{H}.$

\textit{ii)} $Q=-\frac{1}{4}M_{8}^{t}\Gamma \left( Q^{\ast }\right)
M_{8}.\medskip $

\textbf{Proof.} \ i) Let $\ Q$ be a complex quaternion. From Proposition 2.1
i) and ii), we obtain:\newline
$\Gamma ^{t}\left( Q^{\ast }\right) M_{8}=\left( 
\begin{array}{cc}
\lambda \left( x^{\ast }\right) & \lambda \left( y\right) \\ 
-\lambda \left( y^{\ast }\right) & \lambda \left( x\right)%
\end{array}%
\right) \left( 
\begin{array}{c}
\theta M \\ 
-M%
\end{array}%
\right) =\medskip $\newline
=$\left( 
\begin{array}{c}
\lambda \left( x^{\ast }\right) \theta M-\lambda \left( y\right) M \\ 
-\lambda \left( y^{\ast }\right) \theta M-\lambda \left( x\right) M%
\end{array}%
\right) $=$\left( 
\begin{array}{c}
\lambda \left( x^{\ast }i-y\right) M \\ 
-\lambda \left( y^{\ast }i+x\right) M%
\end{array}%
\right) =\medskip $\newline
=$\left( 
\begin{array}{c}
\lambda \left( ix+iiy\right) M \\ 
-\lambda \left( iy+x\right) M%
\end{array}%
\right) $=$\left( 
\begin{array}{c}
\lambda (i\left( x+iy\right) )M \\ 
-M(x+iy)%
\end{array}%
\right) =\medskip $\newline
=$\left( 
\begin{array}{c}
\theta \lambda (x+iy)M \\ 
-M(x+iy)%
\end{array}%
\right) $=$\left( 
\begin{array}{c}
\theta M\left( x+iy\right) \\ 
-M(x+iy)%
\end{array}%
\right) \left( 
\begin{array}{c}
\theta M \\ 
-M%
\end{array}%
\right) \left( x+iy\right) $=$M_{8}Q.\medskip $

ii) If we multiply the relation $\Gamma ^{t}\left( Q^{\ast }\right) M_{8}$=$%
M_{8}Q$ to the left \ side with\medskip \newline
$-\frac{1}{4}M_{8}^{t},$ we obtain $Q$=$-\frac{1}{4}M_{8}^{t}\Gamma
^{t}\left( Q^{\ast }\right) M_{8}.\Box \medskip $

\textbf{Proposition 2.8.} \textit{For} $X,A\in \mathbb{H}_{C},X$=$x+iy,A$=$%
a+ib,$ \textit{we have }%
\begin{equation*}
\Theta \left( XA\right) =\Theta \left( A\right) \Theta (X).
\end{equation*}%
$\medskip $

\textbf{Proof.} Using Proposition 1.1 ii), Proposition 2.2 iv), relations
1.3 and\newline
1.4, it results that \newline
$\Theta \left( XA\right) $=$\left( 
\begin{array}{cc}
\rho \left( xa\text{-}y^{\ast }b\right) & \text{-}\rho \left( x^{\ast }b%
\text{+}ya\right) \\ 
\rho \left( \left( x^{\ast }b\text{+}ya\right) ^{\ast }\right) & \rho \left(
\left( xa\text{-}y^{\ast }b\right) ^{\ast }\right)%
\end{array}%
\right) $=$\medskip \medskip $\newline
=$\left( 
\begin{array}{cc}
\rho \left( xa\text{-}y^{\ast }b\right) & \text{-}\rho \left( x^{\ast }b%
\text{+}ya\right) \\ 
\rho \left( \left( x^{\ast }b\text{+}ya\right) ^{\ast }\right) & \rho \left(
\left( xa\text{-}y^{\ast }b\right) ^{\ast }\right)%
\end{array}%
\right) $=$\medskip \medskip $\newline
=$\left( 
\begin{array}{cc}
\rho \left( xa\text{-}y^{\ast }b\right) & \text{-}\rho \left( x^{\ast }b%
\text{+}ya\right) \\ 
\rho \left( xb^{\ast }\text{+}y^{\ast }a^{\ast }\right) & \rho \left(
x^{\ast }a^{\ast }\text{-}yb^{\ast }\right)%
\end{array}%
\right) .\medskip \medskip $\newline
$\Theta \left( A\right) \Theta \left( X\right) $=$\left( 
\begin{array}{cc}
\rho \left( a\right) & \text{-}\rho \left( b\right) \\ 
\rho \left( b^{\ast }\right) & \rho \left( a^{\ast }\right)%
\end{array}%
\right) \left( 
\begin{array}{cc}
\rho \left( x\right) & \text{-}\rho \left( y\right) \\ 
\rho \left( y^{\ast }\right) & \rho \left( x^{\ast }\right)%
\end{array}%
\right) $=$\medskip \medskip $\newline
=$\left( 
\begin{array}{cc}
\rho \left( a\right) \rho \left( x\right) \text{-}\rho \left( b\right) \rho
\left( y^{\ast }\right) & \text{-}\rho \left( a\right) \rho \left( y\right) 
\text{-}\rho \left( b\right) \rho \left( x^{\ast }\right) \\ 
\rho \left( b^{\ast }\right) \rho \left( x\right) \text{+}\rho \left(
a^{\ast }\right) \rho \left( y^{\ast }\right) & \text{-}\rho \left( b^{\ast
}\right) \rho \left( y\right) \text{+}\rho \left( a^{\ast }\right) \rho
\left( x^{\ast }\right)%
\end{array}%
\right) $=\newline
\medskip

=$\left( 
\begin{array}{cc}
\rho \left( xa-y^{\ast }b\right) & -\rho \left( x^{\ast }b+ya\right) \\ 
\rho \left( xb^{\ast }+y^{\ast }a^{\ast }\right) & \rho \left( x^{\ast
}a^{\ast }-yb^{\ast }\right)%
\end{array}%
\right) .\Box \medskip$ \medskip

\textbf{Proposition 2.9.} \textit{Let} $X\in \mathbb{H}_{C},X=x+iy,x,y\in 
\mathbb{H},~$\textit{then:}

\textit{i)} $\overrightarrow{X}=\left( 
\begin{array}{cc}
1 & 0 \\ 
0 & \alpha%
\end{array}%
\right) \Theta \left( X\right) \left( 
\begin{array}{c}
1 \\ 
0%
\end{array}%
\right) ,$ \textit{where} $1=I_{4}\in \mathcal{M}_{4}\left( \mathbb{R}%
\right) $ \textit{is the identity matrix, } $0=O_{4}\in \mathcal{M}%
_{4}\left( \mathbb{R}\right) $ \textit{is the zero matrix and }$\alpha \in 
\mathcal{M}_{4}\left( \mathbb{R}\right) $ \textit{as in Proposition 2.5 iii)}%
.\medskip

\textit{ii)} $\overrightarrow{XA}=\left( 
\begin{array}{cc}
1 & 0 \\ 
0 & \alpha%
\end{array}%
\right) \Theta \left( A\right) \left( 
\begin{array}{cc}
1 & 0 \\ 
0 & \alpha%
\end{array}%
\right) \overrightarrow{X}.\medskip $

\textit{iii)} $\Gamma \left( A\right) \left( 
\begin{array}{cc}
1 & 0 \\ 
0 & \alpha%
\end{array}%
\right) \Theta \left( B\right) \left( 
\begin{array}{cc}
1 & 0 \\ 
0 & \alpha%
\end{array}%
\right) $=\newline
=$\left( 
\begin{array}{cc}
1 & 0 \\ 
0 & \alpha%
\end{array}%
\right) \Theta \left( B\right) \left( 
\begin{array}{cc}
1 & 0 \\ 
0 & \alpha%
\end{array}%
\right) \Gamma \left( A\right) ,$ \textit{for all }$A,B\in \mathbb{H}%
_{C}.\medskip $

\textbf{Proof.} i) \ We have $\left( 
\begin{array}{cc}
1 & 0 \\ 
0 & \alpha%
\end{array}%
\right) \Theta \left( X\right) \left( 
\begin{array}{c}
1 \\ 
0%
\end{array}%
\right) =\medskip $\newline
=$\left( 
\begin{array}{cc}
1 & 0 \\ 
0 & \alpha%
\end{array}%
\right) \left( 
\begin{array}{cc}
\rho \left( x\right) & -\rho \left( y\right) \\ 
\rho \left( y^{\ast }\right) & \rho \left( x^{\ast }\right)%
\end{array}%
\right) \left( 
\begin{array}{c}
1 \\ 
0%
\end{array}%
\right) $=\medskip \newline
=$\left( 
\begin{array}{cc}
1 & 0 \\ 
0 & \alpha%
\end{array}%
\right) \left( 
\begin{array}{c}
\rho \left( x\right) \\ 
\rho \left( y^{\ast }\right)%
\end{array}%
\right) $=$\left( 
\begin{array}{cc}
1 & 0 \\ 
0 & \alpha%
\end{array}%
\right) \left( 
\begin{array}{c}
\overrightarrow{x} \\ 
\overrightarrow{y^{+}}%
\end{array}%
\right) $=\medskip \newline
=$\left( 
\begin{array}{c}
\overrightarrow{x} \\ 
\alpha \overrightarrow{y^{+}}%
\end{array}%
\right) $=$\left( 
\begin{array}{c}
\overrightarrow{x} \\ 
\overrightarrow{y}%
\end{array}%
\right) .\medskip $

ii)$\ \overrightarrow{XA}$= $\left( 
\begin{array}{cc}
1 & 0 \\ 
0 & \alpha%
\end{array}%
\right) \Theta \left( XA\right) \left( 
\begin{array}{c}
1 \\ 
0%
\end{array}%
\right) $=\medskip \newline
=$\left( 
\begin{array}{cc}
1 & 0 \\ 
0 & \alpha%
\end{array}%
\right) \Theta \left( A\right) \Theta \left( X\right) \left( 
\begin{array}{c}
1 \\ 
0%
\end{array}%
\right) $=\medskip \newline
=$\left( 
\begin{array}{cc}
1 & 0 \\ 
0 & \alpha%
\end{array}%
\right) \Theta \left( A\right) \left( 
\begin{array}{cc}
1 & 0 \\ 
0 & \alpha%
\end{array}%
\right) \left( 
\begin{array}{cc}
1 & 0 \\ 
0 & \alpha%
\end{array}%
\right) \Theta \left( X\right) \left( 
\begin{array}{c}
1 \\ 
0%
\end{array}%
\right) $=\medskip \newline
=$\left( 
\begin{array}{cc}
1 & 0 \\ 
0 & \alpha%
\end{array}%
\right) \Theta \left( A\right) \left( 
\begin{array}{cc}
1 & 0 \\ 
0 & \alpha%
\end{array}%
\right) \overrightarrow{X}.\medskip $

iii) We obtain $\overrightarrow{AXB}$=$\overrightarrow{A(XB)}$=$\Gamma
\left( A\right) \overrightarrow{XB}$=\medskip \newline
=$\Gamma \left( A\right) \left( 
\begin{array}{cc}
1 & 0 \\ 
0 & \alpha%
\end{array}%
\right) \Theta \left( B\right) \left( 
\begin{array}{cc}
1 & 0 \\ 
0 & \alpha%
\end{array}%
\right) \overrightarrow{X}.\medskip $ \newline
Since \ $\overrightarrow{AXB}=\overrightarrow{A(XB)}=\overrightarrow{(AX)B},$
it results that\medskip\ \newline
$\overrightarrow{(AX)B}=\left( 
\begin{array}{cc}
1 & 0 \\ 
0 & \alpha%
\end{array}%
\right) \Theta \left( B\right) \left( 
\begin{array}{cc}
1 & 0 \\ 
0 & \alpha%
\end{array}%
\right) \overrightarrow{AX}$=\medskip \newline
=$\left( 
\begin{array}{cc}
1 & 0 \\ 
0 & \alpha%
\end{array}%
\right) \Theta \left( B\right) \left( 
\begin{array}{cc}
1 & 0 \\ 
0 & \alpha%
\end{array}%
\right) \Gamma \left( A\right) \overrightarrow{X},$ therefore we obtain the
asked relation. $\Box \medskip $

\textbf{Theorem 2.10. }\textit{With the above notations, the following
relation is true:}

\begin{equation*}
\Gamma ^{t}\left( X\right) =M_{1}\Theta \left( X\right) M_{2},
\end{equation*}%
\textit{where}\newline
$M_{1}$=$\left( 
\begin{array}{cc}
-A_{1} & 0 \\ 
0 & A_{1}%
\end{array}%
\right) \in \mathcal{M}_{8}\left( \mathbb{R}\right) ,\medskip ~$\newline
$M_{2}$=$\left( 
\begin{array}{cc}
-A_{2} & 0 \\ 
0 & A_{2}%
\end{array}%
\right) \in \mathcal{M}_{8}\left( \mathbb{R}\right) $ \textit{and\medskip } 
\newline
$A_{1}$=$\left( 
\begin{array}{cccc}
0 & -1 & 0 & 0 \\ 
-1 & 0 & 0 & 0 \\ 
0 & 0 & 0 & 1 \\ 
0 & 0 & -1 & 0%
\end{array}%
\right) \in \mathcal{M}_{4}\left( \mathbb{R}\right) ,\medskip $\newline
$A_{2}$=$\left( 
\begin{array}{cccc}
0 & -1 & 0 & 0 \\ 
-1 & 0 & 0 & 0 \\ 
0 & 0 & 0 & -1 \\ 
0 & 0 & 1 & 0%
\end{array}%
\right) \in \mathcal{M}_{4}\left( \mathbb{R}\right) .\medskip \medskip $

\textbf{Proof.} First, we remark that $A_{1}\rho \left( a\right)
A_{2}=\lambda ^{t}\left( a\right) .$ Indeed,\medskip \newline
$\left( 
\begin{array}{cccc}
0 & \text{-}1 & 0 & 0 \\ 
\text{-}1 & 0 & 0 & 0 \\ 
0 & 0 & 0 & 1 \\ 
0 & 0 & \text{-}1 & 0%
\end{array}%
\right) \left( 
\begin{array}{llll}
a_{0} & \text{-}a_{1} & \text{-}a_{2} & \text{-}a_{3} \\ 
a_{1} & a_{0} & a_{3} & \text{-}a_{2} \\ 
a_{2} & \text{-}a_{3} & a_{0} & a_{1} \\ 
a_{3} & a_{2} & \text{-}a_{1} & a_{0}%
\end{array}%
\right) \left( 
\begin{array}{cccc}
0 & \text{-}1 & 0 & 0 \\ 
\text{-}1 & 0 & 0 & 0 \\ 
0 & 0 & 0 & \text{-}1 \\ 
0 & 0 & 1 & 0%
\end{array}%
\right) $=\medskip \newline
$\left( 
\begin{array}{llll}
-a_{1} & -a_{0} & -a_{3} & a_{2} \\ 
-a_{0} & a_{1} & a_{2} & a_{3} \\ 
a_{3} & a_{2} & -a_{1} & a_{0} \\ 
-a_{2} & a_{3} & -a_{0} & -a_{1}%
\end{array}%
\right) \left( 
\begin{array}{cccc}
0 & -1 & 0 & 0 \\ 
-1 & 0 & 0 & 0 \\ 
0 & 0 & 0 & -1 \\ 
0 & 0 & 1 & 0%
\end{array}%
\right) $=\medskip \newline
$=\left( 
\begin{array}{llll}
a_{0} & a_{1} & a_{2} & a_{3} \\ 
-a_{1} & a_{0} & a_{3} & -a_{2} \\ 
-a_{2} & -a_{3} & a_{0} & a_{1} \\ 
-a_{3} & a_{2} & -a_{1} & a_{0}%
\end{array}%
\right) =\lambda ^{t}\left( a\right) .\medskip $

We have\medskip \newline
$M_{1}\Theta \left( \overline{X}\right) M_{2}=$\newline
$=\left( 
\begin{array}{cc}
-A_{1} & 0 \\ 
0 & A_{1}%
\end{array}%
\right) \left( 
\begin{array}{cc}
\rho \left( x\right) & -\rho \left( y\right) \\ 
\rho \left( y^{\ast }\right) & \rho \left( x^{\ast }\right)%
\end{array}%
\right) \left( 
\begin{array}{cc}
-A_{2} & 0 \\ 
0 & A_{2}%
\end{array}%
\right) =\medskip $\newline
=$\left( 
\begin{array}{cc}
-A_{1}\rho \left( x\right) & A_{1}\rho \left( y\right) \\ 
A_{1}\rho \left( y^{\ast }\right) & A_{1}\rho \left( x^{\ast }\right)%
\end{array}%
\right) \left( 
\begin{array}{cc}
-A_{2} & 0 \\ 
0 & A_{2}%
\end{array}%
\right) =\medskip $\newline
=$\left( 
\begin{array}{cc}
A_{1}\rho \left( x\right) A_{2} & A_{1}\rho \left( y\right) A_{2} \\ 
-A_{1}\rho \left( y^{\ast }\right) A_{2} & A_{1}\rho \left( x^{\ast }\right)
A_{2}%
\end{array}%
\right) $=$\left( 
\begin{array}{cc}
\lambda \left( x\right) & \lambda \left( y\right) \\ 
-\lambda \left( y^{\ast }\right) & \lambda \left( x^{\ast }\right)%
\end{array}%
\right) $=\medskip \newline
=$\left( 
\begin{array}{cc}
\lambda \left( x\right) & -\lambda \left( y^{\ast }\right) \\ 
\lambda \left( y\right) & \lambda \left( x^{\ast }\right)%
\end{array}%
\right) ^{t}=\Gamma ^{t}\left( x\right) .\Box \medskip \medskip $

\textbf{Remark 2.11}. From Theorem 2.7 and Theorem 2.10, it results that

\begin{equation*}
Q=-\frac{1}{4}N_{1}\Theta ^{t}\left( X^{\ast }\right) N_{2},
\end{equation*}%
where $Q\in \mathbb{H}_{C}$ is \textit{\ }a complex quaternion, $N_{1}=$ $%
M_{8}^{t}M_{2}^{t}$ and $N_{2}=M_{1}^{t}M_{8}.\medskip $

\textbf{Proposition 2.12.} \textit{For} $Q\in \mathbb{H}_{C},Q=a+ib,$ 
\textit{we have: }%
\begin{equation*}
\det \Gamma \left( Q\right) \text{=}\det \Theta \left( Q\right) \text{=}%
n\left( aa^{\ast }+b^{\ast }b\right) ^{2}\text{=}n\left( a^{\ast }a+b^{\ast
}b\right) ^{2}.
\end{equation*}

$\medskip $

\textbf{Proof. }

We obtain: $\det \Gamma \left( Q\right) $=$\det \left( 
\begin{array}{cc}
\lambda \left( a\right) & -\lambda \left( b^{\ast }\right) \\ 
\lambda \left( b\right) & \lambda \left( a^{\ast }\right)%
\end{array}%
\right) $=\newline
=$\det \left( \lambda \left( a\right) \lambda \left( a^{\ast }\right) \text{+%
}\lambda \left( b^{\ast }\right) \lambda \left( b\right) \right) $=\newline
=$\det \left( \lambda \left( aa^{\ast }\text{+}b^{\ast }b\right) \right)
=n\left( aa^{\ast }\text{+}b^{\ast }b\right) ^{2}.$

For the second, we have: $\det \Theta \left( Q\right) $=$\det \left( 
\begin{array}{cc}
\rho \left( a\right) & -\rho \left( b\right) \\ 
\rho \left( b^{\ast }\right) & \rho \left( a^{\ast }\right)%
\end{array}%
\right) $=\newline
=$\det \left( \rho \left( a\right) \rho \left( a^{\ast }\right) \text{+}\rho
\left( b\right) \rho \left( b^{\ast }\right) \right) $=\newline
=$\det \left( \rho \left( a^{\ast }a\text{+}b^{\ast }b\right) \right) $=$%
n\left( a^{\ast }a\text{+}b^{\ast }b\right) ^{2}.$

By straightforward calculation, it results that \ $n\left( aa^{\ast }\text{+}%
b^{\ast }b\right) ^{2}$=\newline
=$n\left( a^{\ast }a\text{+}b^{\ast }b\right) ^{2}$. 
\begin{equation*}
\end{equation*}

\textbf{3. Examples}

\begin{equation*}
\end{equation*}

The following sequence of numbers%
\begin{equation*}
0,1,1,2,3,5,8,13,21,....,
\end{equation*}%
with the $n$th term given by the formula:%
\begin{equation*}
f_{n}=f_{n-1}+f_{n-2,}\ n\geq 2,\ 
\end{equation*}%
where $f_{0}=0,f_{1}=1,$ is called the \textit{Fibonacci numbers}$.$

In [Ho; 63], the author defined \ and studied Fibonacci quaternions given by
the formula:%
\begin{equation*}
F_{n}=f_{n}\cdot 1+f_{n+1}e_{2}+f_{n+2}e_{3}+f_{n+3}e_{4},
\end{equation*}%
where $f_{n}$ is the Fibonacci numbers,

\begin{equation*}
e_{m}^{2}=-1,\,\,\,m\in \{2,3,4\}
\end{equation*}%
and \ 
\begin{equation*}
e_{m}e_{q}=-e_{q}e_{m}=\beta _{mq}e_{t},\,\,\beta _{mq}\in \{-1,1\},m\neq
q,m,q\in \{\,\,2,3,4\},
\end{equation*}%
$\ \beta _{mq}$ and $e_{t}$ being uniquely determined by $e_{m}$ and $e_{q}.$
$F_{n}$ is called the \ $n$th Fibonacci quaternion. In the same paper, the
author gave \ some relations for the \ $n$th Fibonacci quaternions, as for
example the norm formula:

\begin{equation*}
n\left( F_{n}\right) =F_{n}\overline{F}_{n}=3f_{2n+3},
\end{equation*}%
where \ $\overline{F}_{n}=f_{n}\cdot
1-f_{n+1}e_{2}-f_{n+2}e_{3}-f_{n+3}e_{4} $ is the conjugate of the $F_{n}.$

In the same paper, Horadam defined the \ $n$th complex Fibonacci numbers as
follows:

\begin{equation*}
q_{n}=f_{n}+if_{n+1},i^{2}=-1,
\end{equation*}%
where $f_{n}$ is the \ $n$th Fibonacci number.

Similarly, the $n$th complex Fibonacci quaternion is the element%
\begin{equation*}
Q_{n}=F_{n}+iF_{n+1},i^{2}=-1,
\end{equation*}%
where \ $F_{n}$ is the \ $n$th Fibonacci quaternion.\medskip

\textbf{Example 3.1. }$\ \ $For the real Fibonacci quaternion $F_{n},$ we
have 
\begin{equation*}
\det \left( \lambda \left( F_{n}\right) \right) =\det \left( \rho \left(
F_{n}\right) \right) =\left( n\left( F_{n}\right) \right) ^{2}=9f_{2n+3}^{2}.
\end{equation*}

\medskip

\textbf{Example 3.2.} The left matrix representation for a complex Fibonacci
quaternion is the matrix: 
\begin{equation*}
\Gamma \left( Q_{n}\right) \text{=}\left( 
\begin{array}{cccccccc}
f_{n} & -f_{n+1} & -f_{n+2} & -f_{n+3} & -f_{n+1} & f_{n+2} & \text{-}f_{n+3}
& \text{-}f_{n+4} \\ 
f_{n+1} & f_{n} & -f_{n+3} & f_{n+2} & -f_{n+2} & \text{-}f_{n+1} & \text{-}%
f_{n+4} & f_{n+3} \\ 
f_{n+2} & f_{n+3} & f_{n} & -f_{n+1} & f_{n+3} & f_{n+4} & \text{-}f_{n+1} & 
f_{n+2} \\ 
f_{n+3} & -f_{n+2} & f_{n+1} & f_{n} & f_{n+4} & \text{-}f_{n+3} & \text{-}%
f_{n+2} & \text{-}f_{n+1} \\ 
f_{n+1} & -f_{n+2} & -f_{n+3} & -f_{n+4} & f_{n} & \text{-}f_{n+1} & f_{n+2}
& f_{n+3} \\ 
f_{n+2} & f_{n+1} & -f_{n+4} & f_{n+3} & f_{n+1} & f_{n} & f_{n+3} & \text{-}%
f_{n+2} \\ 
f_{n+3} & f_{n+4} & f_{n+1} & -f_{n+2} & -f_{n+2} & \text{-}f_{n+3} & f_{n}
& \text{-}f_{n+1} \\ 
f_{n+4} & -f_{n+3} & f_{n+2} & f_{n+1} & -f_{n+3} & f_{n+2} & f_{n+1} & f_{n}%
\end{array}%
\right) .
\end{equation*}%
By straightforward calculation, the determinant of the matrix $\Gamma \left(
Q_{n}\right) $ is\newline
$\det \Gamma \left( Q_{n}\right) $=$\allowbreak \left( f_{n}^{2}\text{+}%
2f_{n}f_{n+2}\text{+}2f_{n+2}^{2}\text{+}f_{n+4}^{2}\text{+}%
2f_{n+2}f_{n+4}\right) ^{2}\allowbreak \cdot $\newline
$\cdot \left( f_{n}^{2}\text{-}2f_{n}f_{n+2}\text{+}4f_{n+1}^{2}\text{+}%
2f_{n+2}^{2}\text{+}4f_{n+3}^{2}\text{+}f_{n+4}^{2}\text{-}%
2f_{n+2}f_{n+4}\right) ^{2}\allowbreak $=\newline
=$\left( (f_{n}\text{+}f_{n+2})^{2}\text{+}\left( f_{n+2}\text{+}%
f_{n+4}\right) ^{2}\right) ^{2}\cdot $\newline
$\cdot \left( \left( f_{n+2}\text{-}f_{n}\right) ^{2}\text{+}\left( f_{n+4}%
\text{-}f_{n+2}\right) ^{2}\text{+}4f_{n+1}^{2}\text{+}4f_{n+3}^{2}\right)
^{2}$=\newline
$=\left( (f_{n}+f_{n+2})^{2}\text{+}\left( f_{n+2}\text{+}f_{n+4}\right)
^{2}\right) ^{2}\left( 5f_{n+1}^{2}\text{+}5f_{n+3}^{2}\right) ^{2}$=\newline
=$25\left( (f_{n}+f_{n+2})^{2}+\left( f_{n+2}+f_{n+4}\right) ^{2}\right)
^{2}\left( f_{n+1}^{2}+f_{n+3}^{2}\right) ^{2}.\medskip $

\textbf{Example 3.3. }The right matrix representation for a complex
Fibonacci quaternion is the matrix:

\begin{equation*}
\Theta \left( Q_{n}\right) \text{=}\left( 
\begin{array}{cccccccc}
f_{n} & \text{-}f_{n+1} & \text{-}f_{n+2} & \text{-}f_{n+3} & \text{-}f_{n+1}
& f_{n+2} & f_{n+3} & f_{n+4} \\ 
f_{n+1} & f_{n} & f_{n+3} & \text{-}f_{n+2} & \text{-}f_{n+2} & \text{-}%
f_{n+1} & \text{-}f_{n+4} & f_{n+3} \\ 
f_{n+2} & \text{-}f_{n+3} & f_{n} & f_{n+1} & \text{-}f_{n+3} & f_{n+4} & 
\text{-}f_{n+1} & \text{-}f_{n+2} \\ 
f_{n+3} & f_{n+2} & \text{-}f_{n+1} & f_{n} & \text{-}f_{n+4} & \text{-}%
f_{n+3} & f_{n+2} & \text{-}f_{n+1} \\ 
f_{n+1} & \text{-}f_{n+2} & f_{n+3} & f_{n+4} & f_{n} & \text{-}f_{n+1} & 
f_{n+2} & f_{n+3} \\ 
f_{n+2} & f_{n+1} & \text{-}f_{n+4} & f_{n+3} & f_{n+1} & f_{n} & \text{-}%
f_{n+3} & f_{n+2} \\ 
\text{-}f_{n+3} & f_{n+4} & f_{n+1} & f_{n+2} & \text{-}f_{n+2} & f_{n+3} & 
f_{n} & f_{n+1} \\ 
\text{-}f_{n+4} & \text{-}f_{n+3} & \text{-}f_{n+2} & f_{n+1} & \text{-}%
f_{n+3} & \text{-}f_{n+2} & \text{-}f_{n+1} & f_{n}%
\end{array}%
\right) .
\end{equation*}%
We have $\ \det \Gamma \left( Q_{n}\right) $=$\allowbreak \left(
f_{n}^{2}+2f_{n}f_{n+2}\text{+}2f_{n+2}^{2}+f_{n+4}^{2}\text{+}%
2f_{n+2}f_{n+4}\right) ^{2}\allowbreak \cdot $\newline
$\cdot \left( f_{n}^{2}\text{-}2f_{n}f_{n+2}\text{+}4f_{n+1}^{2}\text{+}%
2f_{n+2}^{2}\text{+}4f_{n+3}^{2}\text{+}f_{n+4}^{2}\text{-}%
2f_{n+2}f_{n+4}\right) ^{2}\allowbreak $=\newline
=$25\left( (f_{n}+f_{n+2})^{2}+\left( f_{n+2}+f_{n+4}\right) ^{2}\right)
^{2}\left( f_{n+1}^{2}+f_{n+3}^{2}\right) ^{2}.$\medskip $\allowbreak $

\textbf{Remark 3.4.} A matrix representation for the complex Fibonacci
quaternion was introduced \ in [Ha; 12]. This matrix representation, denoted
in the following with $\varepsilon ,$ is a pseudo-representation since $%
\varepsilon \left( XA\right) \neq \varepsilon \left( X\right) \varepsilon
\left( A\right) $ or $\varepsilon \left( XA\right) \neq \varepsilon \left(
A\right) \varepsilon \left( X\right) ,$ where $X,A\in \mathbb{H}%
_{C},X=x+iy,A=a+ib.$ Indeed, using the above notations, we can write the
representation from [Ha; 12] under the form%
\begin{equation*}
\varepsilon \left( A\right) =\left( 
\begin{array}{cc}
\rho ^{t}\left( a\right) & \rho ^{t}\left( b\right) \\ 
-\rho ^{t}\left( b\right) & \rho ^{t}\left( a\right)%
\end{array}%
\right) .
\end{equation*}%
By straightforward calculation, we have%
\begin{equation*}
\varepsilon \left( XA\right) =\left( 
\begin{array}{cc}
\rho ^{t}\left( xa-y^{\ast }b\right) & \rho ^{t}\left( x^{\ast }b+ya\right)
\\ 
-\rho ^{t}\left( x^{\ast }b+ya\right) & \rho ^{t}\left( xa-y^{\ast }b\right)%
\end{array}%
\right) ,
\end{equation*}%
\begin{equation*}
\varepsilon \left( X\right) \varepsilon \left( A\right) =\left( 
\begin{array}{cc}
\rho ^{t}\left( xa-yb\right) & \rho ^{t}\left( xb+ya\right) \\ 
-\rho ^{t}\left( xb+ya\right) & \rho ^{t}\left( xa-yb\right)%
\end{array}%
\right) ,
\end{equation*}%
and 
\begin{equation*}
\varepsilon \left( A\right) \varepsilon \left( X\right) =\left( 
\begin{array}{cc}
\rho ^{t}\left( ax-by\right) & \rho ^{t}\left( bx+ay\right) \\ 
-\rho ^{t}\left( bx+ay\right) & \rho ^{t}\left( ax-by\right)%
\end{array}%
\right) .
\end{equation*}%
\medskip \medskip

From Fundamental Theorem of Algebra, it is known that any polynomial of
degree $n$ with coefficients in a field $\ K$ has at most $n$ roots in $K$.
If the coefficients are in $\mathbb{H}$ (the division real quaternion
algebra), the situation is different. For $\mathbb{H}$ over the real field,
there it is a kind of a fundamental theorem of algebra: \textit{If a
polynomial has only one term of the greatest degree in }$\mathbb{H}$\textit{%
\ then it has at least one root in} $\mathbb{H}$. (see [Ei, Ni; 44] and [Sm;
04]). \ 

In the following, we will give two examples of complex quaternion equations
with more than one greatest term with a unique solution or without solutions.

\medskip

\textbf{Example 3.5.} Let $Q_{n}=F_{n}+iF_{n+1}$ be a complex Fibonacci
quaternion and $A$ a complex quaternion. We consider equations:

\begin{equation}
Q_{n}X-XQ_{n}=A  \tag{3.1.}
\end{equation}%
and 
\begin{equation}
Q_{n}X+XQ_{n}=A  \tag{3.2.}
\end{equation}

If the equation $\left( 3.1\right) $ has a solution, then this solution is
not unique, but the equation $\left( 3.2\right) $ has a unique solution.
Indeed, using the vector representation, Proposition 2.5 and Proposition
2.9, equation $\left( 3.1\right) $ becomes:%
\begin{equation*}
\left( \Gamma \left( Q_{n}\right) -\left( 
\begin{array}{cc}
1 & 0 \\ 
0 & \alpha 
\end{array}%
\right) \Theta \left( Q_{n}\right) \left( 
\begin{array}{cc}
1 & 0 \\ 
0 & \alpha 
\end{array}%
\right) \right) \overrightarrow{X}=\overrightarrow{A}.
\end{equation*}%
We obtain that the matrix $B=$ $\Gamma \left( Q_{n}\right) -\left( 
\begin{array}{cc}
1 & 0 \\ 
0 & \alpha 
\end{array}%
\right) \Theta \left( Q_{n}\right) \left( 
\begin{array}{cc}
1 & 0 \\ 
0 & \alpha 
\end{array}%
\right) $=\newline
=$\left( 
\begin{array}{cccccccc}
0 & 0 & 0 & 0 & 0 & 0 & 0 & 0 \\ 
0 & 0 & \text{-}2f_{n+3} & 2f_{n+2} & 0 & 0 & \text{-}2f_{n+4} & 2f_{n+3} \\ 
0 & 2f_{n+3} & 0 & \text{-}2f_{n+1} & 2f_{n+3} & 0 & \text{-}2f_{n+1} & 0 \\ 
0 & \text{-}2f_{n+2} & 2f_{n+1} & 0 & 2f_{n+4} & 0 & 0 & \text{-}2f_{n+1} \\ 
0 & 0 & \text{-}2f_{n+3} & \text{-}2f_{n+4} & 0 & 0 & 2f_{n+2} & 2f_{n+3} \\ 
0 & 0 & 0 & 0 & 0 & 0 & 0 & 0 \\ 
0 & 2f_{n+4} & 2f_{n+1} & 0 & \text{-}2f_{n+2} & 0 & 0 & \text{-}2f_{n+1} \\ 
0 & \text{-}2f_{n+3} & 0 & 2f_{n+1} & \text{-}2f_{n+3} & 0 & 2f_{n+1} & 0%
\end{array}%
\right) \medskip $ \newline
has $\det B=0$ and $rankB=4,$ as we can find by straightforward calculation.
Therefore, if the equation $\left( 3.1\right) $ has a solution, this
solution is not unique.

In the same way, the equation $\left( 3.2\right) $ becomes 
\begin{equation*}
\left( \Gamma \left( Q_{n}\right) +\left( 
\begin{array}{cc}
1 & 0 \\ 
0 & \alpha%
\end{array}%
\right) \Theta \left( Q_{n}\right) \left( 
\begin{array}{cc}
1 & 0 \\ 
0 & \alpha%
\end{array}%
\right) \right) \overrightarrow{X}=\overrightarrow{A}.
\end{equation*}%
We obtain that the matrix $D=$ $\Gamma \left( Q_{n}\right) +\left( 
\begin{array}{cc}
1 & 0 \\ 
0 & \alpha%
\end{array}%
\right) \Theta \left( Q_{n}\right) \left( 
\begin{array}{cc}
1 & 0 \\ 
0 & \alpha%
\end{array}%
\right) =$\newline
=$\left( 
\begin{array}{cccccccc}
2f_{n} & \text{-}2f_{n+1} & \text{-}2f_{n+2} & \text{-}2f_{n+3} & \text{-}%
2f_{n+1} & 2f_{n+2} & \text{-}2f_{n+3} & \text{-}2f_{n+4} \\ 
2f_{n\text{+}1} & 2f_{n} & 0 & 0 & \text{-}2f_{n+2} & \text{-}2f_{n+1} & 0 & 
0 \\ 
2f_{n\text{+}2} & 0 & 2f_{n} & 0 & 0 & 2f_{n+4} & 0 & 2f_{n+2} \\ 
2f_{n\text{+}3} & 0 & 0 & 2f_{n} & 0 & \text{-}2f_{n+3} & \text{-}2f_{n+2} & 
0 \\ 
2f_{n\text{+}1} & \text{-}2f_{n+2} & 0 & 0 & 2f_{n} & \text{-}2f_{n+1} & 0 & 
0 \\ 
2f_{n\text{+}2} & 2f_{n+1} & \text{-}2f_{n+4} & 2f_{n+3} & 2f_{n+1} & 2f_{n}
& 2f_{n+3} & \text{-}2f_{n+2} \\ 
2f_{n\text{+}3} & 0 & 0 & \text{-}2f_{n+2} & 0 & \text{-}2f_{n+3} & 2f_{n} & 
0 \\ 
2f_{n\text{+}4} & 0 & 2f_{n+2} & 0 & 0 & 2f_{n+2} & 0 & 2f_{n}%
\end{array}%
\right) \medskip $ \newline
has $\ \det D=$\newline
=$256\left( f_{n}\text{-}f_{n+2}\right) ^{2}\left( f_{n}\text{+}%
f_{n+2}\right) ^{2}\left( f_{n}^{2}\text{+}2f_{n}f_{n+2}\text{+}2f_{n+2}^{2}%
\text{+}f_{n+4}^{2}\text{+}2f_{n\text{+}2}f_{n\text{+}4}\right) \cdot $%
\newline
$\cdot \allowbreak \left( f_{n}^{2}\text{-}2f_{n}f_{n+2}\text{+}4f_{n+1}^{2}%
\text{+}2f_{n+2}^{2}\text{+}4f_{n+3}^{2}\text{+}f_{n+4}^{2}\text{-}%
2f_{n+2}f_{n+4}\right) \allowbreak $=\newline
=$1280f_{n+1}^{2}\left( f_{n}\text{+}f_{n+2}\right) ^{2}\left( \left( f_{n}%
\text{+}f_{n+2}\right) ^{2}\text{+}\left( f_{n+2}\text{+}f_{n+4}\right)
^{2}\right) \left( f_{n+1}^{2}\text{+}f_{n+3}^{2}\right) .$ \newline
It results $\det D\neq 0,$ therefore the equation $\left( 3.2\right) $ has a
unique solution.\medskip

\textbf{Example 3.6. }With the above notations, the matrix$\ $%
\begin{equation*}
\ \delta \left( Q_{n}\right) =\Gamma \left( Q_{n}\right) -\Theta \left(
Q_{n}\right)
\end{equation*}%
is an invertible matrix. \newline
Indeed, $\delta \left( Q_{n}\right) =$ \newline
=$\left( 
\begin{array}{cccccccc}
0 & 0 & 0 & 0 & 0 & 0 & \text{-}2f_{n+3} & \text{-}2f_{n+4} \\ 
0 & 0 & \text{-}2f_{n+3} & 2f_{n+2} & 0 & 0 & 0 & 0 \\ 
0 & 2f_{n+3} & 0 & \text{-}2f_{n+1} & 2f_{n+3} & 0 & 0 & 2f_{n+2} \\ 
0 & \text{-}2f_{n+2} & 2f_{n+1} & 0 & 2f_{n+4} & 0 & \text{-}2f_{n+2} & 0 \\ 
0 & 0 & \text{-}2f_{n+3} & \text{-}2f_{n+4} & 0 & 0 & 0 & 0 \\ 
0 & 0 & 0 & 0 & 0 & 0 & 2f_{n+3} & \text{-}2f_{n+2} \\ 
2f_{n+3} & 0 & 0 & \text{-}2f_{n+2} & 0 & \text{-}2f_{n+3} & 0 & \text{-}%
2f_{n+1} \\ 
2f_{n+4} & 0 & 2f_{n+2} & 0 & 0 & 2f_{n+2} & 2f_{n+1} & 0%
\end{array}%
\right) $ \medskip \newline
and 
\begin{equation*}
\det \delta \left( Q_{n}\right) =256\left( f_{n+3}\right) ^{4}\left(
f_{n+2}+f_{n+4}\right) ^{4}
\end{equation*}%
is different from zero.\medskip

\textbf{Conclusions.} In this paper we introduced two real matrix
representation for the complex quaternions and we investigated some of the
properties of these representations. Because of their various applications
to complex quaternions and to \ matrices of complex quaternions, this paper
can be regarded as a starting point for a further research  of
these representations. 
\begin{equation*}
\end{equation*}%
\textbf{Acknowledgements.} Authors thank referee for his/her suggestions
which help us to improve this paper.%
\begin{equation*}
\end{equation*}

\textbf{References}%
\begin{equation*}
\end{equation*}

[Ei, Ni; 44] S.Eilenberg, I.Niven, \textit{The \textquotedblleft\
fundamental theorem of algebra\textquotedblright\ for quaternions}, Bull.
Amer. Math. Soc., \textbf{50}(1944), 246-248.

[Ha; 12] S. Halici, \textit{On complex Fibonacci Quaternions}, Adv. in Appl.
Clifford Algebras, DOI 10.1007/s00006-012-0337-5.

[Ho; 63] A. F. Horadam, \textit{Complex Fibonacci Numbers and Fibonacci
Quaternions}, Amer. Math. Monthly, \textbf{70}(1963), 289-291.

[Sm; 04] W.D.Smith, \textit{Quaternions, octonions, and now, 16-ons, and
2n-ons; New \ kinds of numbers}, \newline
www. math. temple.edu/wds/homepage/nce2.ps, 2004.

[Ti; 00] Y. Tian, \textit{Matrix reprezentations of octonions and their
applications, }Adv. in Appl. Clifford Algebras, \textbf{10}(1)( 2000), 61-90.

[Ti; 00(1)] Y. Tian, \textit{Matrix Theory over the Complex Quaternion
Algebra}, arXiv:math/0004005v1, 1 April 2000.%
\begin{equation*}
\end{equation*}%
\qquad\ \ \ \ \ \ 

Cristina FLAUT

{\small Faculty of Mathematics and Computer Science,}

{\small Ovidius University,}

{\small Bd. Mamaia 124, 900527, CONSTANTA,}

{\small ROMANIA}

{\small http://cristinaflaut.wikispaces.com/}

{\small http://www.univ-ovidius.ro/math/}

{\small e-mail:}

{\small cflaut@univ-ovidius.ro}

{\small cristina\_flaut@yahoo.com}%
\begin{equation*}
\end{equation*}

Vitalii \ SHPAKIVSKYI

{\small Department of Complex Analysis and Potential Theory}

{\small \ Institute of Mathematics of the National Academy of Sciences of
Ukraine,}

{\small \ 3, Tereshchenkivs'ka st.}

{\small \ 01601 Kiev-4}

{\small \ UKRAINE}

{\small \ http://www.imath.kiev.ua/\symbol{126}complex/}

{\small \ e-mail: shpakivskyi@mail.ru}

\end{document}